\documentclass[a4paper,12pt]{article}
\usepackage{amsfonts,amssymb}
\usepackage{graphics}
\title{A numerical scheme for singular shock solutions and a study of its consistence in the sense of distributions}
\author{M. Colombeau,\\ \ \ \\
D\'epartement de Math\'ematiques,\\ Laboratoire Ceregmia,\\Universit\'e des Antilles et de la Guyane, \\97157, Pointe-\`a-Pitre Cedex, Guadeloupe, France\\ \& \\ Instituto de Matematica e Estastistica,\\Universidade de Sao Paulo,\\Cidade Universitaria, SP, Brazil}

\begin{document}
\maketitle

\begin {abstract} 

In this paper we present a numerical scheme for the approximation of singular shock solutions of the Keyfitz-Kranzer model system. Consistence in the sense of distributions is studied.  As long as some numerical properties are verified when the space step tends to 0, we prove that the scheme provides a numerical solution that satisfies the equations in the sense of distributions with an approximation that tends to 0 when $h\rightarrow 0$. We also show that this scheme adapts to degenerate systems. This is illustrated by two examples:  the system presenting delta wave solutions originally studied by Korchinski and another system studied by Keyfitz-Kranzer that models elasticity. Consistence of the scheme in the sense of distributions is fully proved in the case of the Korchinski  model.

\end{abstract}

\textbf{1. Introduction}. Singular shocks have been discovered  and investigated by different authors, see in particular [1,4,5,6]. They have been observed from various viscosity techniques: Dafermos-Di Perna viscosity in  [1], usual viscosity in [4]. In the case of singular shocks viscosity solutions converge so weakly that their pointwise limits do not satisfy the classical Rankine-Hugoniot conditions. Besides this fact a unique entropic solution to the Riemann problem has been obtained in [1] for arbitrarily large data.  In this paper we propose a numerical scheme based on a splitting technique that captures the singular shocks. We observe   results exactly similar to those obtained in [1,4] with their respective viscosity techniques. Studies have shown the relevance of this scheme for other systems presenting irregular solutions. In our study of irregular shocks we consider two standard first order model systems of two equations whose solutions of the Riemann problem involve singular shocks and delta shocks.  We also notice that this scheme provides neat results for the  Keyfitz-Kranzer system of elasticity [2] for which the intrinsic difficulty is different from those in the two  systems above.\\

 This paper focusses on the Keyfitz-Kranzer system 
 
\begin{equation}   u_t+(u^2-v)_x=0,    \end{equation}
\begin{equation} v_t+(\frac{1}{3} u^3-u)_x=0,       \end{equation}
which produces singular shocks [1,4], and   the  system 
\begin{equation}   u_t+(u^2)_x=0,    \end{equation}
\begin{equation} v_t+(uv)_x=0,       \end{equation}
originally considered by Korchinski [3] who  discovered and investigated delta shocks in the solution of the Riemann problem.\\ 

Let $u_h,v_h$ be the sequence of approximate solutions from the scheme. Under simple numerical properties to be rigorously proved, or  to be admitted from numerical tests, we prove that the scheme 
 is consistent in the sense of distributions in the following sense: $\forall (\phi,\psi) \in (\mathcal{C}_c ^\infty(\mathbb{R}\times \mathbb{R}^+))^2$,
   
  \begin{equation}  \int[u_h \phi_t+((u_h)^2-v_h)\phi_x] dxdt \rightarrow 0, \ \ \int[v_h \psi_t+(\frac{1}{3} (u_h)^3-u_h)\psi_x] dxdt \rightarrow 0,\end{equation} \\
  respectively  $\int[u_h \phi_t+((u_h)^2)\phi_x] dxdt \rightarrow 0, \ \ \int[v_h \psi_t+(u_hv_h)\psi_x] dxdt \rightarrow 0,$\\
  \\
when the space step $h\rightarrow 0$. This means that the functions $u_h,v_h$ tend to satisfy the equations when $h\rightarrow 0$.\\

 For system (1,2) we  check numerically that the needed assumptions are satisfied for values of $h$ as small as possible.  We rigorously prove that, in the case of system (3,4), for any initial condition  $u^0\in L^1(\mathbb{R})\cap L^\infty (\mathbb{R})$ and $v^0\in L^1(\mathbb{R})$, these assumptions are satisfied. Therefore the scheme is consistent in the above sense. Of course, in the first case,  from a rigorous point of view, one cannot be sure that these numerical assumptions always hold for every $h$ when $h\rightarrow 0$. The proof in this paper shows that, for any given family of test functions with uniformly bounded support and uniformly bounded first and second derivatives, then a numerical solution satisfies the equations in the sense of distributions within a small deviation  depending on $h$ whenever these assumptions remain valid. \\

\textbf{2. A  numerical scheme}.   The singular shocks of the Keyfitz-Kranzer equations are unbounded which makes  the elaboration of numerical schemes difficult: in the scheme below the numerical velocity $u$ in system (1,2)  can be unbounded when the space step $h$ tends to 0 which forces us to accept that the CFL coefficient $r$ tends to 0 when $h\rightarrow 0$ in order to preserve the CFL condition $r\|u\|_{L^\infty}\leq 1$. Therefore $r=r_h$ depends on $h$ and also on time so that $r_h\|u_h\|_{L^\infty}\leq 1$.\\ 

If $r_h$ tends to 0 ( i.e. if  $\|u_h\|_{L^\infty}$ tends to $\infty$) slowly enough, then one can nevertheless obtain a convenient numerical scheme, although of an order less than one, on condition that for each iteration the assumptions are verified when $h\rightarrow 0$. This ensures consistence of the scheme in the sense of distributions, although the limit is not a distribution in general: it can be a singular shock in the case of the Keyfitz-Kranzer equations.   Numerical results are given to prove that   the set of assumptions  is satisfied in representative situations of singular shocks.
 In the case of the Keyfitz-Kranzer equations the scheme consists in a splitting  of equations into the two subsystems 

\begin{equation} u_t +( u^2)_x=0,      \end{equation}
\begin{equation}  v_t+(v u)_x=0,     \end{equation}
\\
which is treated by transport with velocity $u$,  and

\begin{equation} u_t =v_x,      \end{equation}
\begin{equation}  v_t=(vu-\frac{u^3}{3}+u)_x,     \end{equation}
\\
which is treated  by a centered discretization. In between, we introduce an average step  in $u,v$ which is needed in general to avoid oscillations due to the centered discretization. More generally the method applies to  systems

\begin{equation} u_t +[u\Phi(u,v)]_x=[A(u,v)]_x,      \end{equation}
\begin{equation} v_t +[v\Phi(u,v)]_x=[B(u,v)]_x,      \end{equation}
\\
which are split into the two subsystems

\begin{equation} u_t +[u\Phi(u,v)]_x=0,      \end{equation}
\begin{equation} v_t +[v\Phi(u,v)]_x=0,      \end{equation}
 where $\Phi(u,v)$ plays the role of numerical velocity and 
 \begin{equation} u_t =[A(u,v)]_x,      \end{equation}
\begin{equation} v_t =[B(u,v)]_x.      \end{equation}

 Systems (12,13) is a family of degenerate systems considered in [2]. In particular the scheme in this paper gives neat results for the system (4) in reference [2] which models an elastic string problem.\\ 

\textit{ The numerical scheme.} The real line is divided into intervals $I_i=]ih-\frac{1}{2} h,ih+\frac{1}{2} h[,i\in \mathbb{Z}$. We set $t_n=nrh$ for $r$ small enough. We will construct step functions  $u(x,t)$, $v(x,t)$ depending on $h$, which are constant on the rectangles $I_i\times ]t_n,t_{n+1}[$, whose step values are denoted $ u_i^n, v_i^n$ respectively. The indices $h$ are often skipped to simplify the notation: $u$ stands for $u_h,...$. If $a<b$ one sets  
\begin{equation} L(a,b):= length \ of \ [0,1]\cap[a,b],     \end{equation}
i.e.
\begin{equation}L(a,b)=max(0,min(1,b)-max(0,a)).\end{equation}
The notation $L$ allows a synthetic formulation of the transport, without being forced to distinguish several cases depending on the signs of the numerical velocities.
 By induction we assume that the set of values $\{u_i^n,v_i^n\}_{i\in \mathbb{Z}}$ is known. We obtain the set of values $\{u_i^{n+1},v_i^{n+1}\}_{i\in \mathbb{Z}}$ as follows.\\

$\bullet$\textit{First step: transport with velocity $\Phi$ during time $rh$}
\begin{equation} \Phi_i^n:=\Phi(u_i^n,v_i^n),    \end{equation}
\begin{equation}  \overline{u}_i:=u_{i-1}^nL(-1+r\Phi_{i-1}^n,r\Phi_{i-1}^n)+ u_{i}^nL(r\Phi_{i}^n,1+r\Phi_{i}^n)+u_{i+1}^nL(1+r\Phi_{i+1}^n,2+r\Phi_{i+1}^n),    \end{equation}

\begin{equation}  \overline{v}_i:=v_{i-1}^nL(-1+r\Phi_{i-1}^n,r\Phi_{i-1}^n)+ v_{i}^nL(r\Phi_{i}^n,1+r\Phi_{i}^n)+v_{i+1}^nL(1+r\Phi_{i+1}^n,2+r\Phi_{i+1}^n).    \end{equation}
\\
When the CFL condition $r|\Phi_i ^n|\leq 1 \ \  \forall i, \ \forall n$ is satisfied the first terms in (19,20), when multiplied by $h$,  represent the quantities $u,v$ issued from the cell $I_{i-1}$ between times $t_n$ and $t_{n+1}$ that lie in the cell $I_i$ at time $t_{n+1}$. Indeed the cell $I_{i-1}=[(i-\frac{3}{2})h,(i-\frac{1}{2})h]$ has been transported according to the vector $r\Phi_{i-1}^nh$, since $\Phi_{i-1}^n$ is the numerical velocity and the duration time is $rh$. The overlap with the fixed cell $I_i=[(i-\frac{1}{2})h,(i+\frac{1}{2})h]$ has a length of $r\Phi_{i-1}^nh$ if $\Phi_{i-1}^n\geq 0$, 0 if $\Phi_{i-1}^n\leq 0$, taking into account the CFL condition $r|\Phi_{i-1}^n|\leq 1$. From (16) one finds $L(-1+r\Phi_{i-1}^n,r\Phi_{i-1}^n)= r\Phi_{i-1}^n$ if $\Phi_{i-1}^n\geq 0$, 0 if $\Phi_{i-1}^n\leq 0$. Division by $h$ is due to the fact that $\overline{u}_i, u_j^n$ are mean values on cells of length $h$.\\
\\
 The second terms in (19,20), when multiplied by $h$, represent the quantities $u,v$ issued from the cell $I_i$ that remain in $I_i$ at time $t_{n+1}$. Indeed the cell $[(i-\frac{1}{2})h,(i+\frac{1}{2})h]$ has been transported by the vector $r\Phi_i^n h$. The overlap with the fixed cell $[(i-\frac{1}{2})h,(i+\frac{1}{2})h]$ is $h-r\Phi_i^n h$ if $\Phi_i^n\geq 0$, $h+r\Phi_i^n h$ if $\Phi_i^n\leq 0$. From (16) one finds $L(r\Phi_{i}^n,1+r\Phi_{i}^n)=1-r\Phi_i^n$ if $\Phi_i^n\geq 0$, $1+r\Phi_i^n$ if $\Phi_i^n\leq 0$. \\
 \\
 The third terms are similar to the first ones: they concern the quantities $u,v$  issued from the cell $I_{i+1}$ that lie in the cell $I_i$ at time $t_{n+1}$, with the same verification as above. Note that $\overline{u}_i,\overline{v}_i$ depend on $n$, which is not explicitely stated to shorten the notation. \\

$\bullet$\textit{Averaging step.} For a value $\alpha, 0\leq\alpha<0.5$, to be chosen, we set
   
 \begin{equation} \widetilde{u}_i:=\alpha \overline{u}_{i-1}+(1-2\alpha)\overline{u}_i+\alpha\overline{u}_{i+1},     \end{equation}
 
  \begin{equation} \widetilde{v}_i:=\alpha \overline{v}_{i-1}+(1-2\alpha)\overline{v}_i+\alpha\overline{v}_{i+1}.     \end{equation}
\\
In the case $A=0,B=0$ the averaging step is useless. Indeed the idea underlying the elaboration of the scheme is that the first step works  well without averaging, and that the numerical defects of the centered discretization in the last step should be compensated by the averaging step performed before it. The splitting should be chosen so as to minimize the importance of the terms involved in the last step.\\

 $\bullet$\textit{Last step: centered discretization}

\begin{equation} u_i^{n+1}:= \widetilde{u}_i+\frac{r}{2}[A(u_{i+1}^n,v_{i+1}^n)-A(u_{i-1}^n,v_{i-1}^n)],\end{equation}

\begin{equation} v_i^{n+1}:= \widetilde{v}_i+\frac{r}{2}[B(u_{i+1}^n,v_{i+1}^n)-B(u_{i-1}^n,v_{i-1}^n)].    \end{equation}\\

  The scheme works well for singular shocks and delta shocks. The theorem below shows that it gives an approximate solution of the equations.\\

\textbf{3. Statement of the theorem}.
 Let $T>0$ be given. Let us seek a solution on $\mathbb{R}\times[0,T]$. The initial conditions $u^0,v^0$ are discretized as usual by mean values in the cells since they are supposed to be $L^1$ functions. Let us apply the scheme under the assumptions (25-29) below: there exists a sequence of values $h$, $h\rightarrow 0$,  a corresponding sequence of values  $r$, $r>0$, and real numbers
 $\beta, \gamma \in [0,1[ $  such that when $h\rightarrow 0 $

\begin{equation}\frac{h}{r}\rightarrow 0      \end{equation}
\begin{equation} \forall n\leq \frac{T}{rh} \ \forall i \ \ \ r|\Phi_i^n| \leq 1,      \end{equation}
 which is the CFL condition, 
 \begin{equation}\forall n\leq \frac{T}{rh} \ \forall i \ \ \  h^\beta |\Phi_i^n|=O(1),      \end{equation}
 which is a constraint on the numerical velocity allowing it to tend to infinity,
 
\begin{equation}\forall n\leq \frac{T}{rh} \ \forall i \ \ \ \sum_i|u_i^n|h=O(1), \ \sum_i|v_i^n|h=O(1),      \end{equation}
which is the $L^1$-stability in $u,v$,
\begin{equation} \forall n\leq \frac{T}{rh} \ \forall i \ \ \  \sum_i|A(u_i^n, v_i^n)|h^{1+\gamma}=O(1),\sum_i|B(u_i^n, v_i^n)|h^{1+\gamma}=O(1).   \end{equation}
\\
 \textbf{Theorem. Consistence of the  scheme.} \textit{ As long as (25-29) are satisfied 
then the scheme is consistent on $\mathbb{R}\times ]0,T[$ in the sense of distributions, i.e. if $u_h,v_h$,  are the step functions from the scheme, then}, $\forall \psi \in \mathcal{C}_c^\infty(\mathbb{R}\times ]0,T[)$, 
  
 \begin{equation} \int{[u_h\psi_t+u_h\Phi(u_h,v_h)\psi_x-A(u_h,v_h)\psi_x]dxdt} \rightarrow 0,      \end{equation}
\begin{equation} \int{[v_h\psi_t+v_h\Phi(u_h,v_h)\psi_x-B(u_h,v_h)\psi_x]dxdt} \rightarrow 0,      \end{equation}
when $h\rightarrow 0.$ \textit{ More precisely the integrals in (30,31) are equal to}
\begin{equation} O(\frac{h}{r})+O(h^{1-\beta})+O (h^{1-\gamma}).      \end{equation}

The scheme will be of order one in the usual cases in which $r$ is constant, $\beta=\gamma=0$, but of an order strictly less than one for singular shocks from the fact that the values of the numerical velocity increase when $h\rightarrow 0$, which forces  $r\rightarrow 0$ and $\beta>0$.\\

 \textbf{4. Proof of the theorem.}  $\int{u\psi_tdxdt}=\sum_{i,n}u_i^n\int_{cell_{i,n}}{\psi_tdxdt}=\sum_{i,n}u_i^n[(\psi_t)_i^n+O(h)]rh^2=\sum_{i,n}u_i^n\frac{\psi_i^n-\psi_i^{n-1}}{rh}rh^2+\sum_{i,n}u_i^nO(rh)rh^2+\sum_{i,n}u_i^nO(h)rh^2.$\\
 Since $|\sum_{i,n}u_i^nO(h)rh^2|\leq \sum_nrh\sum_i|u_i^n||O(h)|h\leq const. T|O(h)|$ from (28), one obtains\\
\begin{equation}  \int{u\psi_tdxdt}=  \sum_{i,n}(u_i^n-u_i^{n+1})h\psi_i^n+O(h). \end{equation}
\\
Similarly\\
\\
$\int{\Phi(u,v)u\psi_xdxdt}=\sum_{i,n}\Phi_i^nu_i^n\int_{cell_{i,n}}{\psi_x dxdt}=\sum_{i,n}\Phi_i^nu_i^n(\psi_x)_i^nrh^2+\sum_{i,n}\Phi_i^nu_i^nO(h)rh^2=\sum_{i,n}\Phi_i^nu_i^n\frac{\psi_{i+1}^n-\psi_i^n}{h}rh^2+\sum_{i,n}\Phi_i^nu_i^nO(h)rh^2+\sum_{i,n}\Phi_i^nu_i^nO(h)rh^2.$\\

$|\sum_{i,n}\Phi_i^nu_i^nO(h)rh^2|\leq \sum_nrh\sum_i|\Phi_i^n||u_i^n||O(h)|h\leq const.T h^{-\beta}h\leq const.h^{1-\beta}$ from (27,28). 
Finally\\
\begin{equation} \int{\Phi(u,v)u\psi_xdxdt}=-h\sum_{i,n}(\Phi_i^nu_i^n-\Phi_{i-1}^nu_{i-1}^n)r\psi_i^n+O(h^{1-\beta}).\end{equation}
Similarly\\

$\int A(u,v)\psi_xdxdt=\sum_{i,n}A(u_i^n,v_i^n)\int_{cell_{i,n}}\psi_xdxdt=\sum_{i,n}A(u_i^n,v_i^n)(\psi_x)_i^nrh^2+\sum_{i,n}A(u_i^n,v_i^n)O(h)rh^2=\sum_{i,n}A(u_i^n,v_i^n)\frac{\psi_{i+1}^n-\psi_i^n}{h}rh^2+\sum_{i,n}A(u_i^n,v_i^n)O(h)rh^2.$\\

$|\sum_{i,n}A(u_i^n,v_i^n)O(h)rh^2|\leq \sum_n rh \sum_i |A(u_i^n,v_i^n)||O(h)|h\leq const.T.h^{-\gamma}h\leq const.h^{1-\gamma}$ from (29).
Therefore
\begin{equation}\int A(u,v)\psi_xdxdt=\sum_{i,n}rh[A(u_{i-1}^n,v_{i-1}^n)-A(u_i^n,v_i^n)]\psi_i^n+O(h^{1-\gamma}).       \end{equation}
\\
 Setting 
 \begin{equation} I:=\int[u\psi_t+u\Phi(u,v)\psi_x-A(u,v)\psi_x]dxdt,      \end{equation}
 one finally obtains from (33-36)
 $$I=-h\sum_{i,n}[u_i^{n+1}-u_i^n+r(u_i^n\Phi_i^n-u_{i-1}^n\Phi_{i-1}^n)-$$\begin{equation}r(A(u_i^n,v_i^n)-A(u_{i-1}^n,v_{i-1}^n))]\psi_i^n +O(h)+O(h^{1-\beta})+O(h^{1-\gamma}).\end{equation}\\
 Up to this point the formulas of the scheme have not yet been used. From (23,21)
 
 $u_i^{n+1}=\overline{u}_i+\alpha(\overline{u}_{i-1}-2\overline{u}_i+\overline{u}_{i+1})+\frac{r}{2}[A(u_{i+1}^n,v_{i+1}^n)-A(u_{i-1}^n,v_{i-1}^n)] $.\\ 
  Therefore, from (37)
\begin{equation} I=I_1+I_2+I_3+O(h)+O(h^{1-\beta})+O(h^{1-\gamma}),     \end{equation}
where
\begin{equation} I_1=  -h\sum_{i,n}[\overline{u}_i-u_i^n+r(u_i^n\Phi_i^n-u_{i-1}^n\Phi_{i-1}^n)]\psi_i^n,  \end{equation}
\begin{equation} I_2=-h\alpha\sum_{i,n} (\overline{u}_{i-1}-2\overline{u}_i+\overline{u}_{i+1})\psi_i^n,     \end{equation}  
 \begin{equation}  I_3=-\frac{1}{2}\sum_{i,n}hr\{ A(u_{i+1}^n,v_{i+1}^n)-A(u_{i-1}^n,v_{i-1}^n)-2 [A(u_{i}^n,v_{i}^n)-A(u_{i-1}^n,v_{i-1}^n)] \}\psi_i^n.  \end{equation}
 We are going to prove successively bounds for  $I_1,I_2,I_3$.\\
 
 $\bullet$ \textit{ bound for $I_1$.} In $I_1$ fix an index $i_0$ and consider successively the two cases $\Phi_{i_0}^n\leq 0$ and $\Phi_{i_0}^n\geq 0$.\\
 
 If $\Phi_{i_0}^n\leq 0$ then, from (16) and the CFL condition (26), $L(r\Phi_{i_0}^n,1+r\Phi_{i_0}^n)=1+r\Phi_{i_0}^n$, $L(1+r\Phi_{i_0}^n,2+r\Phi_{i_0}^n)=-r\Phi_{i_0}^n$ and $L(-1+r\Phi_{i_0}^n,r\Phi_{i_0}^n)=0$ . Therefore from (19)\\
\\ 
 $\overline{u}_{i_0}=u_{i_0}^n(1+r\Phi_{i_0}^n)$+terms not involving $u_{i_0}^n,$\\
 $\overline{u}_{i_0-1}=-u_{i_0}^nr\Phi_{i_0}^n$+terms not involving $u_{i_0}^n,$\\
 $\overline{u}_{i_0+1}$ does not involve $u_{i_0}^n.$\\
 \\
 From the CFL condition the other terms $\overline{u}_i$ do not involve $u_{i_0}^n$. Therefore in the sum $\sum_i\overline{u}_i\psi_i^n$ the term $u_{i_0}^n$ occurs in (and only in) 
 $$u_{i_0}^n(1+r\Phi_{i_0}^n)\psi_{i_0}^n-u_{i_0}^nr\Phi_{i_0}^n\psi_{i_0-1}^n.$$
 \\
 Consequently in the sum  $\sum_{i}[\overline{u}_i-u_i^n+r(u_i^n\Phi_i^n-u_{i-1}^n\Phi_{i-1}^n)]\psi_i^n$, the term involving $u_{i_0}^n$ is
 
 \begin{equation}u_{i_0}^n(1+r\Phi_{i_0}^n)\psi_{i_0}^n-u_{i_0}^nr\Phi_{i_0}^n\psi_{i_0-1}^n-u_{i_0}^n\psi_{i_0}^n+ru_{i_0}^n\Phi_{i_0}^n\psi_{i_0}^n-ru_{i_0}^n\Phi_{i_0}^n\psi_{i_0+1}^n\end{equation}
 \\
where the first two terms come from  $\overline{u}_{i_0}$ and  $\overline{u}_{i_0-1}$. The sum (42) is equal to $ru_{i_0}^n(\Phi_{i_0}^n)[\psi_{i_0}^n-\psi_{i_0-1}^n+\psi_{i_0}^n-\psi_{i_0+1}^n]=ru_{i_0}^n\Phi_{i_0}^nO(h^2)$ from Taylor's formula applied to $\psi$.\\

 If $\Phi_{i_0}^n\geq 0$ then, an analogous reasoning involving  $\overline{u}_{i_0}$ and $\overline{u}_{i_0+1}$ instead of $\overline{u}_{i_0}$ and $\overline{u}_{i_0-1}$ gives the value 0. Therefore from (39)\\

 $|I_1|\leq h\sum_{i_0,n}u_{i_0}^nr\Phi_{i_0}^nO(h^2)=\sum_nrh \sum_i \Phi_i^nu_i^nhO(h),$
 i.e. from (27,28)
\begin{equation} I_1=O(h^{1- \beta}).     \end{equation}

$\bullet$ \textit{ bound for $I_2$.} From (40) $I_2=-h\alpha\sum_{i,n}\overline{u}_i(\psi_{i+1}^n-2\psi_{i}^n+\psi_{i-1}^n)=\alpha\sum_{n}rh \frac{1}{r}\sum_i\overline{u}_iO(h^2)=\alpha T\frac{h}{r} O(1)$
since one has $\sum_i|\overline{u}_i|h\leq\sum_i|u_i^n|h=O(1)$. Indeed (19) implies the formula
\begin{equation}|\overline{u}_i|\leq  |u_{i-1}^n|L(-1+r\Phi_{i-1}^n,r\Phi_{i-1}^n)+ |u_{i}^n|L(r\Phi_{i}^n,1+r\Phi_{i}^n)+|u_{i+1}^n|L(1+r\Phi_{i+1}^n,2+r\Phi_{i+1}^n).\end{equation}
The definition (16) of $L$  implies $L(-1+a,a)+L(a,1+a)+L(1+a,1+2a)=1$. Therefore from (44) $\sum_i|\overline{u}_i|\leq \sum_i|u_i^n|$.
This implies

\begin{equation} I_2=O(\frac{h}{r}).    \end{equation}

$\bullet$  \textit{ bound for $I_3$.}  \ $I_3=-\frac{hr}{2}\sum_{i,n}\{ A(u_{i}^n,v_{i}^n)\psi_{i-1}^n-A(u_{i}^n,v_{i}^n)\psi_{i+1}^n-2 A(u_{i}^n,v_{i}^n)\psi_i^n +2A(u_{i}^n,v_{i}^n) \psi_{i+1}^n\}=-\frac{1}{2}\sum_{n}rh\sum_iA(u_{i}^n,v_{i}^n)[\psi_{i-1}^n-2\psi_{i}^n+\psi_{i+1}^n]=const. Th^{-\gamma}O(h)$ from Taylor's formula in $\psi$ and (29). Therefore

\begin{equation}I_3=O(h^{1-\gamma}).       \end{equation}

Finally from (38,43,45,46) 
  
 \begin{equation} I=O(h^{1-\beta})+O(h^{1-\gamma})+O(\frac{h}{r}),      \end{equation}
 which ends the proof.$\Box$\\

  \textbf{5. Approximation of the Keyfitz-Kranzer system (1,2).}   We consider successively the three different typical solutions of Riemann problems in figures 8,7,6 in  [4]: singular shock, intermediate overcompressive shock and usual shocks. The numerical solutions obtained from the scheme are identical to those shown in [4] even in  absence of additional viscosity. We consider first the Riemann problem  in figure 8 in [4], which shows a singular shock. The initial data is  $(u_l,v_l,u_r,v_r)=(1.5,0,-2.065426,1.410639)$. We adopt the values $\alpha=0.2, \beta=0.5, \gamma=0.4$.  
One  chooses the value of $r_h$ close to  the maximum value of $r$ that satisfies the CFL condition (26). For simplicity we denote 
$$"(27)":=h^\beta max_{i,n}|u_i^n|,$$
  $$"(28)":= max_{n}(\sum_i|u_i^n|h,\sum_i|v_i^n|h  ),$$ 
  $$ "(29)":= max_{n}(\sum_i|A(u_i^n, v_i^n)|h^{1+\gamma}, \sum_i|B(u_i^n, v_i^n)|h^{1+\gamma}  )$$
  for the values in the assumptions of the theorem. \\ 
  
  In order to check the consistence theorem we present the values of $\frac{h}{r}$ that must tend to 0 from (25), and the values "27","28", "29" that must be bounded. Results of a test for $T=5$ with the interval $[-4,4]$ are given in the table below.\\

\begin{tabular}{|c|c|c|c|c|c|c|}
  
 \hline
 $ h$  & $r$ & $\frac{h}{r}$ & $"(27)"$ &$"(28)"$ &$"(29)"$\\
 
 \hline
 0.0400&0.300&0.1333&0.6289&14.97&3.62\\
 \hline
 0.0200&0.240&0.0833&0.5830&14.97&2.84\\
 \hline
 0.0100&0.170&0.0588&0.5309&14.96&2.26\\
 \hline
 0.0050&0.132&0.0379&0.5271&14.93&1.84\\
 \hline
 0.0025&0.095&0.0263&0.5178&14.90&1.53\\
 \hline
 0.00125&0.065&0.0192&0.5021&14.87&1.29\\
 \hline
 0.00062&0.040&0.0156&0.4326&14.85&1.10\\
 \hline
0.00031&0.025&0.0125&0.4024&14.83&0.96\\
\hline
 \end{tabular}\\   \\

 Now we choose $T=1$ and the interval $[-0.5,0.5]$ in order to reach   
  smaller values of $h$. 
 The values of the parameters are again $\alpha=0.2, \beta=0.5, \gamma=0.4$\\

 \begin{tabular}{|c|c|c|c|c|c|c|}

 \hline
 $ h$  & $r$ & $\frac{h}{r}$ & $"(27)"$ &$"(28)"$ &$"(29)"$\\
 
 \hline
 0.0020&0.18&0.0111&0.2444&1.9232&0.1791\\
 \hline
 0.0010&0.13&0.0077&0.2337&1.9170&0.1480\\
 \hline
 0.0005&0.09&0.0056&0.2225&1.9109&0.1252\\
 \hline
 0.00025&0.06&0.0042&0.2090&1.9054&0.1081\\
 \hline
 0.000125&0.043&0.0029&0.2070&1.8999&0.0973\\
 \hline
 0.0000833&0.035&0.0024&0.2051&1.8972&0.0926\\
\hline
0.0000625&0.030&0.0021&0.2028&1.8955&0.0898\\
\hline
0.0000500&0.026&0.0019&0.1979&1.8944&0.0874\\
\hline
0.0000333&0.021&0.0016&0.1955&1.8923&0.0848\\
\hline
0.0000250&0.019&0.0013&0.2010&1.8907&0.0847\\
\hline
0.0000166&0.015&0.0011&0.1957&1.8891&0.0829\\
\hline
0.0000125&0.012&0.0010&0.1847&1.8884&0.0803\\
\hline

 \end{tabular}\\  \\

 \textbf{insert figure 1}\\
 
 \textit{Figure 1: The numerical solution from the last test in the second table.  One can observe that the scheme  reproduces exactly the aspect of the singular shock in figure 8, in reference [4].} \\

  The  values of $r$  are chosen close to the maximum  values for which the scheme satisfies the CFL condition $r\|u\|_\infty \leq 1$ . One observes that the quantity $\frac{h}{r}\rightarrow 0$ as $\sqrt{h}$ and that the three quantities in the columns "27", "28","29" are bounded (since quantity "27" is proportional to the sup. of $|u|$ it is very sensitive to the chosen value of $r$ close to the sup. of values of $r$ that satisfy the CFL condition). Therefore, since $\beta=0.5, \gamma=0.4$, the scheme is of order 0.5 in $h$ from (32). This is not a good result in general from a numerical viewpoint; however the presence of singular shocks gives a numerical velocity which is of the order $\frac{1}{\sqrt h}$  instead of a constant in the usual situations in which the scheme is always of order 1. We can see also that the bounds in the proof of the theorem are not optimal since one has used a bound involving the factor  $\|\Phi_i^n\|_\infty$  while $\Phi_i^n=u_i^n$ is uniformly bounded independently of time except on the singular shock. Indeed  one can see that the scheme gives acceptable results. On a standard PC top values of the peak in $v$ in the above tests have  reached the value $3700$ for the Riemann problem under consideration  while they have reached values $10^{6}$ in the case of system (3,4). For the system (3,4) we will rigorously prove in section 6  that the scheme is of order 1. The set of results in these two tables   gives a reasonable presumption that the decrease of $\frac{h}{ r}$ and the boundedness of the three  quantities "27","28","29" continue to hold when $h\rightarrow 0$, which would allow the theorem to be applied with confidence. If we only consider  values of $h$ for which (25-29) have been tested then the proof of the theorem gives a bound (depending on the sup norm of the derivatives of order two of $\psi$ and its support) for the integrals in (30,31) according to (32). \\
 
 Then we consider the Riemann problem $(u_l,v_l,u_r,v_r)=(1.5, 0, -1.895644, 1.343466)$ represented in figure 7 in [4]. In this case we choose $\alpha=0.2, \beta=0,\gamma=0$. We obtain the following table:\\

 \begin{tabular}{|c|c|c|c|c|c|c|}
 
 \hline
 $ h$  & $r$ & $\frac{h}{r}$ & $"(27)"$ &$"(28)"$ &$"(29)"$\\
 
 \hline
 0.0050&0.45&0.0111&1.9205&1.6913&1.2831\\
 \hline
 0.0010&0.45&0.0022&1.9205&1.6965&1.2753\\
 \hline
 0.0005&0.45&0.0011&1.9205&1.6972&1.2743\\
 \hline
 0.00025&0.45&0.0006&1.9205&1.6975&1.2739\\
 \hline
 0.000125&0.45&0.0003&1.9205&1.6977&1.2736\\
 \hline
 0.0000625&0.45&0.0001&1.9205&1.6977&1.2735\\
\hline
 \end{tabular}\\ \\

 \textbf{insert figure 2}\\
 
 \textit{Figure 2: The numerical solution from the Riemann problem considered   in the third table ($h=0.04, r=0.45$).  One can observe that the scheme in this paper reproduces exactly the aspect of the limit overcompressive shock in  figure 7 in [4]. An enlargement has been done in the  horizontal direction  to observe the detailed structure of the shock.} \\ 
 
The results are very clear due to the boundedness of $u$ in this case.  There is a very natural presumption that these results continue to hold when $h\rightarrow 0$. One can see that the scheme is of order one in $h$ as this follows from the theorem.\\
 
 For the third Riemann problem,  $(u_l,v_l,u_r,v_r)=(1.5, 0,-1.725862, 1.276293)$ in figure 6 in [4], in which there is no singular shock, the results are very clear, exactly the same as  those in the above table. We have always observed  results as good in the case of bounded numerical velocity.\\
 
  We now present a system for which a full proof of consistence in the sense of distributions has been obtained.\\

\textbf{6. Application to system (3,4). }\\ 
 
 $\bullet$ One considers the $2\times 2$ system (3,4)  which produces delta-waves in the variable $v$, see [3]. Here  $\Phi(u,v)=u, \  A=B=0$. In this case one can choose $\alpha=0$ in (21,22) since the last step (23,24) is absent. Then $u_i^{n+1}=\overline{u_i}, v_i^{n+1}=\overline{v_i}$; the choice $\alpha>0$ works as well with the same proofs. It follows from (44) that $\sum_i |\overline{u}_i| \leq \sum_i|u_i^n|$. Therefore by induction on $n  \ \sum_i|u_i^{n+1}|\leq\sum_i |u_i^0|$. The same  proof applies for $v$. Choosing the initial condition $u^0,v^0$ in $L^1$ this proves (28). To prove (25,26,27) we will prove the maximum principle in the numerical velocity $u$.\\

\textit{ Lemma}.  If $r max_i |u_i^0| \leq \frac{1}{2}$ \textit{then $u$ satisfies the maximum principle}.\\ 
\\ 
\textit{Proof.} Let the index $i$ be fixed. Consider the various possible combinations of signs in the three values  $u_{i-1}^n,u_{i}^n,u_{i+1}^n$. In each case one will check that  
 $$min(u_{i-1}^n,u_{i}^n,u_{i+1}^n) \leq \overline{u}_i = u_i^{n+1} \leq max(u_{i-1}^n,u_{i}^n,u_{i+1}^n)$$
which  proves the maximum principle by induction on $n$. By induction up to order $n$ the condition  $r max_i |u_i^0| \leq \frac{1}{2}$ implies  $r max_i |u_i^n| \leq \frac{1}{2}$. Now we pass to order $n+1$.\\
 \\  
 $\bullet$ case (+,+,+). Formula (19) with $\Phi=u$ gives 
 \begin{equation}\overline{u}_i= u_{i-1}^nr u_{i-1}^n+u_i^n(1-ru_i^n)=u_i^n+r(u_{i-1}^n-u_i^n)(u_{i-1}^n+u_i^n).\end{equation}
 First note that $\overline{u}_i \geq 0$ because $1-ru_i^n\geq 0$ from the property $r max_i |u_i^n| \leq \frac{1}{2}$.
 We consider successively the two cases $u_i^n\geq u_{i-1}^n$ and $u_i^n\leq u_{i-1}^n$. If $u_i^n\geq u_{i-1}^n$ then (48) gives 
 $\overline{u}_i\leq u_{i}^n$.  If $u_i^n\leq u_{i-1}^n$ then 
  $\overline{u}_i-u_{i-1}^n=(u_i^n-u_{i-1}^n)[1-r(u_i^n+u_{i-1}^n)]\leq 0$ since the last factor is $\geq 0$ by induction. We have checked that 
 $$0\leq \overline{u}_i \leq max(u_{i-1}^n,u_i^n).$$
 \\
 $\bullet$ case (+,+,-). Formula (19) gives 
 \begin{equation}\overline{u}_i= u_{i-1}^nr u_{i-1}^n+u_i^n(1-ru_i^n)+u_{i+1}^n(-r {u}_{i+1}^n).\end{equation}
 
 First let us prove that $\overline{u}_i\geq u_{i+1}^n$. The properties $u_{i-1}^n\geq 0,u_{i}^n\geq 0,ru_{i}^n\leq \frac{1}{2} $ imply that 
  
 $\overline{u}_i\geq u_{i+1}^n (-ru_{i+1}^n)\geq u_{i+1}^n $ since $0\leq -ru_{i+1}^n\leq \frac{1}{2} $ and $u_{i+1}^n\leq 0. $\\
 
 Now let us check that $\overline{u}_i\leq max(u_{i-1}^n, u_i^n)$. Formula (49) and $u_{i+1}^n\leq 0$ imply $\overline{u}_i \leq u_{i-1}^nr u_{i-1}^n+u_i^n(1-ru_i^n)$. From this inequality the proof is the same as in the case (+++).\\
 \\
 $\bullet$ case (-,+,+). Formula (19) gives $\overline{u}_i=u_i^n(1-ru_i^n)$ which implies $\overline{u}_i\leq  u_i^n$ since $0\leq ru_i^n\leq \frac{1}{2}$ and, $\overline{u}_i\geq 0$.\\
 \\
 $\bullet$ case (-,+,-). Formula (19) gives 
 $$\overline{u}_i=u_i^n(1-ru_i^n)+u_{i+1}^n(-r {u}_{i+1}^n) =u_i^n+r[-(u_{i+1}^n)^2-(u_i^n)^2]\leq u_i^n.$$ 
 Now $\overline{u}_i-
 u_{i+1}^n=u_i^n-u_{i+1}^n-r[(u_{i+1}^n)^2+(u_i^n)^2].$  Since $u_i^n u_{i+1}^n\leq 0, (u_i^n)^2+(u_{i+1}^n)^2\leq (u_i^n)^2+(u_{i+1}^n)^2-2u_i^nu_{i+1}^n=(u_i^n-u_{i+1}^n)^2$. Therefore $\overline{u}_i- u_{i+1}^n \geq u_i^n-u_{i+1}^n-r(u_{i}^n-u_{i+1}^n)^2 = (u_i^n-u_{i+1}^n)[1-r(u_i^n-u_{i+1}^n)]\geq 0$ since the second factor is positive, which implies $\overline{u}_i\geq u_{i+1}^n. $\\
 
 In the four  cases in which $u_i^n\leq 0$ the verifications are similar.\\
 
 Finally we have proved properties (25-28), with  $r$ independent of $h$, $\beta=0$ and $\gamma=0$ since $A=B=0$. Therefore from the theorem  the scheme converges in the sense of distributions and is of order one in $h$.  It has been checked numerically that its real interpretation is the well known  solution.\\

 \textbf{6. Conclusion.}  We have presented a numerical scheme which captures the singular shock solutions of the Keyfitz-Kranzer model without recourse to a vanishing viscosity method. We have  observed numerically  exactly the same results previously observed by the various authors. The consistence of the scheme for this system has been checked numerically up to very small values of $h$. The theorem states that the approximate solutions from the scheme tend to satisfy the equations in the sense of distributions. This scheme adapts to degenerate systems such as the Korchinski model system and the Keyfitz-Kranzer system of elasticity. In the case of the Korchinski system consistence in the sense of distributions has been fully proved.  \\

 \end{document}